\documentclass[12pt,a4paper]{article}

\usepackage{newtxtext}
\usepackage{newtxmath}

\usepackage[margin=1in]{geometry}
\usepackage{amsmath}
\usepackage{graphicx}
\usepackage{hyperref}
\usepackage{bm}        
\usepackage{authblk}   
\usepackage{cite}

\usepackage{titlesec}
\titleformat{\section}{\large\bfseries}{\thesection}{1em}{}
\titleformat{\subsection}{\normalsize\bfseries}{\thesubsection}{1em}{}
\titleformat{\subsubsection}{\normalsize\itshape}{\thesubsubsection}{1em}{}

\usepackage[labelfont=bf,font=small]{caption}

\hypersetup{
	colorlinks=true,
	linkcolor=blue,
	citecolor=blue,
	urlcolor=blue
}

\usepackage{orcidlink} 
\usepackage{balance}
\usepackage{latexsym}
\usepackage{float}
\usepackage{xcolor}

	\title{An explicit integrator uniform in the true anomaly and exactly preserving all integrals of motion in the three-dimensional Kepler problem}

\author{Jan L. Cie\'sli\'nski\,\orcidlink{0000-0003-1730-0950} \footnote{j.cieslinski@uwb.edu.pl} \ } 

\author{Maciej Jurgielewicz\,\orcidlink{0009-0006-0349-1150} \footnote{m.jurgielewicz@uwb.edu.pl}}

\affil{Faculty of Physics, University of Bialystok, ul.\ Ciolkowskiego 1L, 15-245 Bialystok, Poland}

\date{}

\begin{document}

	\newcommand{\red}{\color{red}}
	
	\newcommand{\rf}[1]{(\ref{#1})}
	\newcommand{\rff}[2]{(\ref{#1}\ref{#2})}
	
	\newcommand{\ba}{\begin{array}}
		\newcommand{\ea}{\end{array}}
	
	\newcommand{\be}{\begin{equation}}
		\newcommand{\ee}{\end{equation}}

	\newcommand{\const}{{\rm const}}
	\newcommand{\ep}{\varepsilon}
	\newcommand{\Cl}{{\cal C}}
	\newcommand{\rr}{{\pmb r}}
	\newcommand{\qq}{{\pmb q}}
	\newcommand{\ph}{\varphi}
	\newcommand{\R}{{\mathbb R}}  
	\newcommand{\C}{{\mathbb C}}
	
	\newcommand{\Q}{{\pmb Q}}
	\newcommand{\PP}{{\pmb P}}

	\newcommand{\e}{{\bf e}}
	
	\newcommand{\m}{\left( \ba{r}}
	\newcommand{\ema}{\ea \right)}
	\newcommand{\mm}{\left( \ba{cc}}
	\newcommand{\miv}{\left( \ba{cccc}}
	
	\newcommand{\scal}[2]{\mbox{$\langle #1 \! \mid #2 \rangle $}}
	\newcommand{\ods}{\par \vspace{0.5cm} \par}
	\newcommand{\dis}{\displaystyle }
	\newcommand{\mc}{\multicolumn}
	
	\newtheorem{prop}{Proposition}
	\newtheorem{Th}{Theorem}
	\newtheorem{lem}{Lemma}
	\newtheorem{rem}{Remark}
	\newtheorem{cor}{Corollary}
	\newtheorem{Def}{Definition}
	\newtheorem{open}{Open problem}
	\newtheorem{ex}{Example}
	\newtheorem{exer}{Exercise}
	
	\newenvironment{Proof}{\par \vspace{2ex} \par
		\noindent \small {\it Proof:}}{\hfill $\Box$ 
		\vspace{2ex} \par }

\maketitle

\begin{abstract}
We develop a numerical scheme for the Kepler problem that preserves exactly all first integrals: angular momentum, total energy, and the Laplace–Runge–Lenz vector.
This property ensures that orbital trajectories retain their precise shape and orientation over long times, avoiding the spurious precession typical of many standard methods.
The scheme uses an adaptive time step derived from a constant angular increment. Analytical considerations and numerical experiments demonstrate that the algorithm combines high accuracy with long-term stability.
\end{abstract}


\section{Introduction}

The Kepler problem, describing the motion of a body under an inverse square
central force and the related two-body problem, remains a cornerstone of celestial mechanics and dynamical systems theory \cite{Cor-2003,MD-1999}. Although the analytical solution for the shape of the Keplerian orbits is well known, the numerical integration of Keplerian trajectories in three dimensions continues to present challenges, especially when long-term stability and accuracy are required \cite{HSD-2025}. Standard numerical methods often suffer from qualitative inaccuracies during
long-term simulations, including gradual energy drift, artificial precession
of the orbital ellipse, and distortion of the orbital shape. These numerical
artifacts can obscure physical insights, especially in high-precision applications
such as long-term planetary motion prediction or spacecraft trajectory
design. 

Geometric numerical integration approach have addressed these issues by designing structure preserving or geometrically inspired numerical schemes that aim to conserve specific invariants and geometric structures of the considered problem  \cite{HLW-2006}. Rich geometric structure of the Kepler problem (see, e.g., \cite{Cor-2003,BGK-2022,MPW-2013,IKMS-2023}) 
makes it an excellent candidate for applying geometric numerical methods, see,e.g.,  \cite{Bre-1999,Cie-2010,Ci-AVF-2014,EE-2017,ME-2025}.
The combination of symplecticity  and time reversibility  is responsible for the excellent (although not exact) long-term stability of modern geometric integrators \cite{HLW-2006,Yos-1990,HS-2005}.   Symplectic integrators are widely used in many applications \cite{WH-1991,For-2006}. They have become famous for their use in long-term integration of the dynamics of the Solar System  \cite{WH-1991}. A key feature of the Wisdom–Holman approach is the separation of the Keplerian motion. Therefore, its numerical integration is a key component in multi-body problems where one body has a dominant mass.

However, to preserve the energy by symplectic integrators, time adaptation needs to be used \cite{EEI-2021,SBPW-2022,YCWYL-2025}, which, in general, is a challenging task  \cite{HS-2005,DL-2025}. To preserve exactly (up to round-off errors) other first integrals, another class of geometric integrators is needed, see, e.g.,  \cite{MQR-1998,Koz-2007,CR-2010}. 

In this Letter we present a discretization of the three-dimensional Kepler problem that exactly preserves all first integrals:   the angular momentum, the total energy, and the Laplace–Runge–Lenz vector.
As a consequence, the method maintains the exact shape and orientation of the orbital ellipse over arbitrarily long time intervals, completely eliminating the usual problem of numerical precession.

The proposed numerical scheme employs an adaptive time step determined by a constant angular increment: the step size is relatively large far from the center and becomes much  smaller near periapsis \cite{Cie-2007}.
In other words, the integration is performed uniformly in the true anomaly (the angle between the radius vector and the direction of periapsis), yielding successive position vectors separated by a constant angular interval. This distinctive feature of our algorithm enables the accurate determination of the body’s angular position along its orbit (the true anomaly) and is related to one of the regularization methods of the Kepler problem \cite{Roa-2017}.  

Both theoretical analysis and numerical experiments demonstrate the excellent accuracy and long-term stability of the proposed algorithm.

\section{Main result: An explicit conservative scheme uniform in the true anomaly}

Motivated by results of the paper \cite{Cie-2007} we propose the following numerical scheme for the 3-dimensional Kepler problem
\be  \label{Kepler}
\frac{d  \pmb{p}}{d t} = - \frac{k \qq}{q^3}  \ , \qquad 
\pmb{p} = m \frac{d \qq}{d t} \ ,
\ee
where $\pmb{q}, \pmb{p} \in \R^3$, $q = |\pmb{q}|$ and $m, k$ are constants.

First, given $\pmb{q}_n$, $\pmb{p}_n$ and $h_n$, and denoting $q_n = |\pmb{q}_n|$, we define auxiliary variables $\pmb{S}_n$ and $\rr_n$:
\be   \label{Sn-rn} 
S_n := \frac{h_n \, \pmb{q}_n \cdot \pmb{p}_n}{ m q_n}  \ , \qquad 
\pmb{r}_n = \pmb{q}_n + \frac{h_n}{2 m} \left( \frac{S_n}{q_n + \sqrt{q_n^2 + S_n^2} } - 1  \right) \pmb{p}_n   \ .
\ee
The proposed algorithm for generating the next step ($\pmb{q}_{n+1}$, $\pmb{p}_{n+1}$, $h_{n+1}$) is given by:
\be  \ba{l} \label{Kep-dis} \displaystyle
\rr_{n+1} =  \rr_n + \frac{h_n \pmb{p}_n }{m} 
\ , \\[3ex]
\displaystyle
\pmb{p}_{n+1} = \pmb{p}_n - \frac{k h_n \rr_{n+1} }{   r_{n+1}^2 r_n \cos\delta} \ , \\[4ex]\dis 
h_{n+1} = \frac{h_n }{\displaystyle  \ \frac{2 r_n \cos 2 \delta}{r_{n+1}} - 1 + \frac{k h_n^2}{  m r_{n+1}^2 r_n \cos\delta}} \ ,  \\[7ex]\dis 
\rr_{n+2} =  \rr_{n+1} + \frac{h_{n+1} \pmb{p}_{n+1} }{m}  \ ,  \\[3ex]\dis 
\pmb{q}_{n+1} = \frac{r_{n+2} \rr_{n+1} +  r_{n+1} \rr_{n+2}}{r_{n+1} + r_{n+2}} \ ,
\ea \ee
where $r_n := |\rr_n|$ and $\delta$ is a constant which  can be defined by initial data ($\rr_0$, $\pmb{p}_0$ and $h_0$) as follows:
\be  \label{delta} 
\cos 2 \delta = \frac{r_0^2 + \rr_0 \cdot \pmb{P}_0 }{r_0 \sqrt{r_0^2 + 2   \rr_0 \cdot \pmb{P}_0  +   \pmb{P}_0^2  } } \ , \qquad \pmb{P}_0 := \frac{h_0 \pmb{p}_0}{m} \ .
\ee
In the next section we describe and discuss properties of this algorithm (derivations and proofs are contained in the Appendix).

\section{Properties of the proposed numerical scheme}

The formula \rf{delta} results directly from $\rr_0 \cdot \rr_1 = r_0 r_1 \cos 2\delta$ which means that the angle between $\rr_0$ and $\rr_1$ equals $2\delta$. In order to avoid the vanishing of $\cos 2\delta$ it suffices to assume that $|\pmb{P}_0| < r_0$, which holds for   $h_0$  sufficiently small. This assumption is reasonable and not too restrictive -- very large time steps are of no practical use anyway. 

Although formulas \rf{Sn-rn} are needed only at the first step of the algorithm \rf{Kep-dis} but they are consistent with it, i.e., they are satisfied for  $n$ replaced by $n+1$, compare  appendix~\ref{app-first}. 

From the formula for $\pmb{q}_{n+1}$ given in \rf{Kep-dis} it follows that $\pmb{q}_n$ (for any $n \in \mathbb N$) is the angle–bisector vector between $\rr_n$ and $\rr_{n+1}$ with its tip lying on the segment joining the endpoints of $\rr_n$ and $\rr_{n+1}$. Therefore, the angle between $\pmb{q}_n$ and $\pmb{q}_{n+1}$ is $2 \delta$, as well. In particular we have
\be 
\pmb{q_n} \cdot \rr_n = q_n r_n \cos \delta \ , \qquad
\pmb{q_n} \cdot \rr_{n+1} = q_n r_{n+1} \cos \delta \ .
\ee

It turns  out that the numerical scheme \rf{Kep-dis} is uniform in terms of the true anomaly, i.e, the angles between $\rr_n$, $\pmb{q}_n$ and $\rr_{n+1}$ do not depend on $n$ (for a proof see appendix~\ref{app-true}). What is more, the numerical scheme \rf{Kep-dis} has the same first integrals as the continuum system \rf{Kepler}, namely:   
\begin{itemize}
	\item  Angular momentum 
	\be \label{ang}
	\pmb{L}_n = \pmb{q}_n \times \pmb{p}_n  \ .
	\ee
	\item Energy  
	\be  \label{ene}
	E_n = \frac{(\pmb{p}_n)^2}{2 m} - \frac{k}{q_n}  
	\ .
	\ee
	\item  Laplace-Runge-Lenz (LRL) vector
	\be   \label{LRL} 
	\pmb{A}_n = \frac{\pmb{p}_n \times \pmb{L}_n }{m} - \frac{k\,  \pmb{q}_n}{q_n}  
	\ee
	(the LRL vector is sometimes defined as  $m \pmb{A}_n$ rather than $\pmb{A}_n$). 
\end{itemize}
As shown in the Appendix, $\pmb{L}_n, E_n$ and $\pmb{A}_n$  do not depend on $n$ and can be expressed by initial data $\pmb{q}_0$ and $\pmb{p}_0$.  Therefore they are exactly preserved by the  algorithm \rf{Kep-dis} (up to round-off errors).

Another remarkable property of the scheme \rf{Kep-dis} is that it allows the exact (up to round-off) determination of the epoch $t_n$ corresponding to the position  $\pmb{q}_n$. 
Indeed, our algorithm yields the true anomaly linear in $\delta$
\be \label{nu_n}
\nu_n = 2 n \delta + \nu_0 \ ,
\ee
where $\nu_0$ is the angle between $\pmb{q}_0$ and $\pmb{A}_0$, which points in the direction of the periaxis. Then, we use well known formulas expressing the time $t_n$ in terms of $\nu_n$. 

In the elliptic case ($E_0 <0$)  
the true anomaly $\nu_n$ can be expressed in terms of the eccentric anomaly $u_n$ as follows
\be
\tan \frac{u_n}{2} = \sqrt{\frac{1-e}{1+e}} \ \tan \frac{\nu_n}{2} \ .
\ee
Here $e$ is the orbital eccentricity, which can be expressed by the first integrals as follows
\be
e = \sqrt{ 1 + \frac{2 E_0 L_0^2}{m k^2}} ,
\ee 
where $L_0 = |\pmb{L}_0|$. Then, we  use the celebrated Kepler equation 
\be  \label{Kep eq}
M_n = u_n - e \sin u_n 
\ee
where  $u_n$ is the eccentric anomaly, i.e., the angle measured at the center of the ellipse between the orbit's periapsis and the current position, and $M_n$ is the mean anomaly, which is a linear function of time
\be
M_n =  M_0 + \frac{2 \sqrt{2} |E_0|^{3/2} } {k\sqrt{m}} (t_n -t_0) \ .
\ee
Analogous formulas for the hyperbolic and parabolic cases can be found, for instance, in \cite{Cor-2003} (where natural units, 
$k=m=1$, are used, but physical constants can be restored by dimensional analysis).

In other algorithms for the Kepler problem, anomalies $\nu_n$ and $u_n$ are not known and have to computed  from the time $t_n$ via the Kepler equation \rf{Kep eq}, which is not explicit with respect to $u_n$. Solvers and approaches to the Kepler equation are constantly being improved to address this key issue \cite{ZBL-2022,WZLZ-2023}.

\section{Results of numerical simulations}

We will compare five numerical schemes, including three standard methods and the algorithm proposed in this Letter:  
\begin{itemize}
	\item MTPI (modified trajectory-preserving integrator), given by formulas \rf{Kep-dis},
	\item TPI (trajectory-preserving integrator), see  \cite{Cie-2007}, 
	\item RK 4 : Runge-Kutta method of the fourth order (see, for instance, {\tt run\_kut4} in \cite{Kiu-2013}),
	\item LF: leapfrog method (see \cite{HLW-2006}, section I.1.4),
	\item S-Y(4): Suzuki-Yoshida method of the fourth order (\cite{Suz-1990,Yos-1990}, see also \cite{CR-2011}), 
\end{itemize}

The efficiency of the algorithms will be assessed in terms of their ability to preserve the first integrals and the shape and orientation of the orbit in space.

In numerical simulations we take $k=3$, $m=0.5$ and the following  initial conditions:
\be  \label{inic}  
\pmb{q}(0)=(100,0.0,0.1)  , \quad \pmb{p}(0)=(0,0.01,0) .
\ee
The time step was taken as $h=0.01$ (for LF), $h=0.02$ (for RK4 and S-Y(4)), and $h_0 = 10$ for MTPI and TPI (which corresponds to $\delta = 0.001$). 

For elliptic orbits, the number of integration steps per period  is given by
\begin{equation}\label{N12}
	N_1=\frac{T}{h} \ , \quad 
	N_2=\frac{\pi}{\delta} \ ,
\end{equation}
corresponding to constant time-step and constant angular-increment methods, respectively.

In our simulations the number of steps per period is $4.53\cdot 10^{4}$ for RK4 and S-Y(4), $9.07\cdot 10^{4}$ for LF, and only $3.14 \cdot 10^3$ for TPI and MTPI.

\subsection{Energy conservation}

The energy $E_0$ is computed from \rf{ene} using initial data $\pmb{q}_0$ and $\pmb{p}_0$. The relative error in energy  is computed as follows:  
\be
E_{err} (n) =\sup_{0 \leqslant j \leqslant n} \left|\frac{E_j-E_0}{E_0}\right| \ .
\ee
Note that, because the error is computed over the entire interval, the resulting error function is non-decreasing.
\begin{figure}[H]
	\centering
	\includegraphics[width=0.47\textwidth]{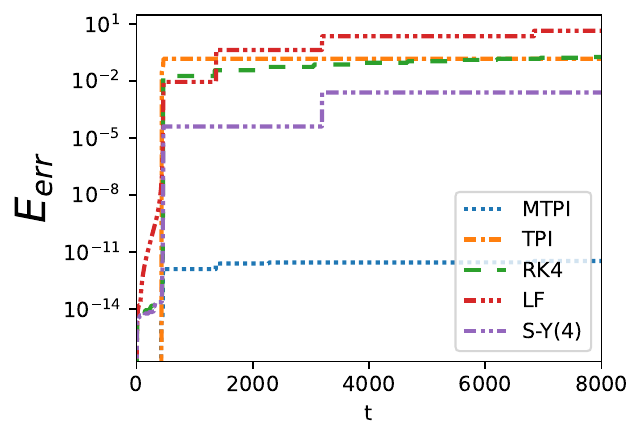} 
	\includegraphics[width=0.47\textwidth]{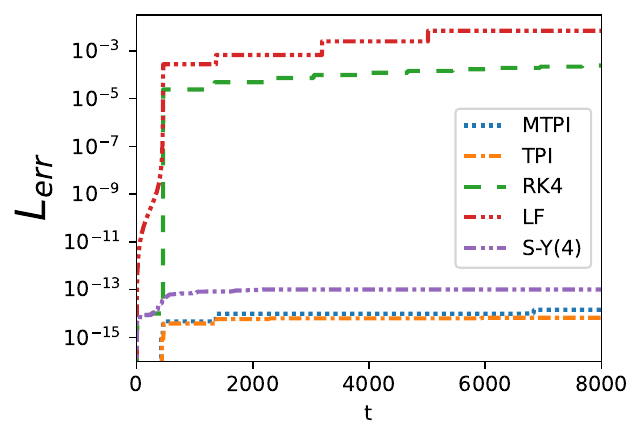}
	\caption{Left: Relative error in energy over time. Right: Relative error in the angular momentum magnitude. }
	\label{fig:E}
\end{figure}

\noindent MTPI shows superior accuracy (by several orders of magnitude)  over all other schemes considered, see Fig.~\ref{fig:E}.

\subsection{Angular momentum conservation}

To quantify the deviation of the angular momentum direction, we use the direction cosine error:
\be \label{dirL}
dir L_{err} (n) =\sup_{0 \leqslant j \leqslant n}
\left(1-\frac{\pmb{L}_j\cdot \pmb{L}_0}{L_j L_0}\right) .
\ee
All the numerical schemes considered pass this test with excellent precision. The deviations remain below $2.3 \cdot 10^{-16}$ 
over the entire interval, confirming that the methods are well suited for orbital motion simulations.

To test the conservation of the angular momentum we compute the relative error in its magnitude:  
\be  \label{Lerr}
L_{err} (n) = \sup_{0 \leqslant j \leqslant n} \left|\frac{L_j-L_0}{L_0}\right|
\ee
Figure~\ref{fig:E} shows excellent performance of TPI and MTPI. Although S-Y(4) performs comparably, its accuracy is nevertheless lower by about one order of magnitude.

\subsection{Laplace-Runge-Lenz vector conservation}

In full analogy to the angular momentum case, we test both the preservation of the vector magnitude and the angular deviation: 
\be  \label{Aerr}
A_{err} (n) =\sup_{0 \leqslant j \leqslant n} \left|\frac{A_j-A_0}{A_0} \right| ,
\ee
\be  \label{dirA}
dir A_{err} (n) =\sup_{0 \leqslant j \leqslant n}
\left(1-\frac{\pmb{A}\cdot \pmb{A}_0}{A_j A_0}\right) .
\ee
\begin{figure}[H]
	\centering
	\begin{tabular}{cc}
		\includegraphics[width=0.47\textwidth]{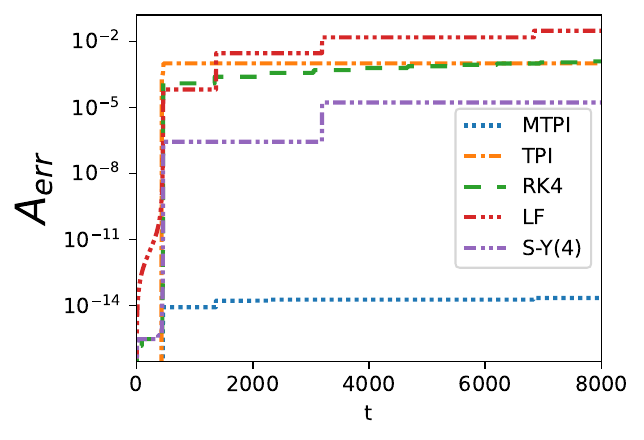} &  \includegraphics[width=0.47\textwidth]{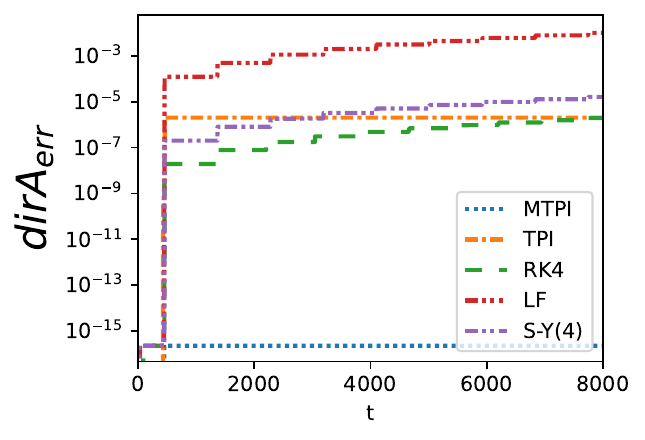} 
	\end{tabular}
	\caption{Left: Relative error in the magnitude of the LRL vector (see \rf{Aerr}). Right: LRL vector direction error (see \rf{dirA}).}
	\label{fig:A}
\end{figure}
In both case MTPI is the best by several orders of magnitude as compared to the other algorithms, see Fig.~\ref{fig:A}.

\subsection{Trajectory preservation}

\begin{figure}[H]
	\centering
	\includegraphics[width=0.5\textwidth]{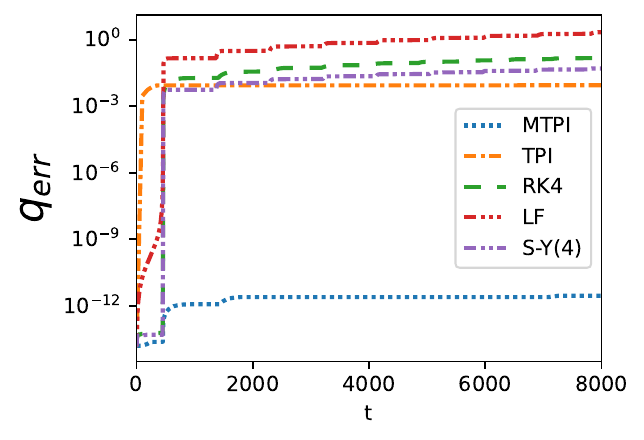} 
	\caption{Relative radial distance error over time. }
	\label{fig:orbite}
\end{figure}

The analytical equation for the orbit shape (see, e.g., \cite{Cor-2003}) allows to compute the distance from the center provided that the true anomaly is given: 
\be \label{orbit}
\frac{1}{q (t_n)}=\frac{k m}{L_0^2}+\frac{k m}{L_0^2}\sqrt{1+\frac{2 E_0 L_0^2}{k^2 m}}\cos(\theta_n + \nu_0) ,
\ee
where $\nu_0$  and $\theta_n$  are given by: 
\be
\cos{\theta_n}=\frac{\pmb{q}_0 \pmb{q}_n}{q_0 q_n} \ , \qquad 
\cos\nu_0 = \frac{\pmb{q}_0 \pmb{A}_0}{q_0 A_0} .
\ee
Thus, $\theta_n$ is obtained from the simulation, while $\nu_0$ comes from the initial data. We define the relative radial distance error as follows:
\be
q_{err} (n) =\sup_{0 \leqslant j \leqslant n} \frac{| q (t_j)-q_j|}{| q (t_j)|} \ .
\ee
Again, MTPI is superior to the other numerical schemes tested, see Fig.~\ref{fig:orbite}. It should be noted that TPI performs better than the three standard algorithms.

\begin{figure}[H]
	\centering
	\includegraphics[width=0.95\textwidth]{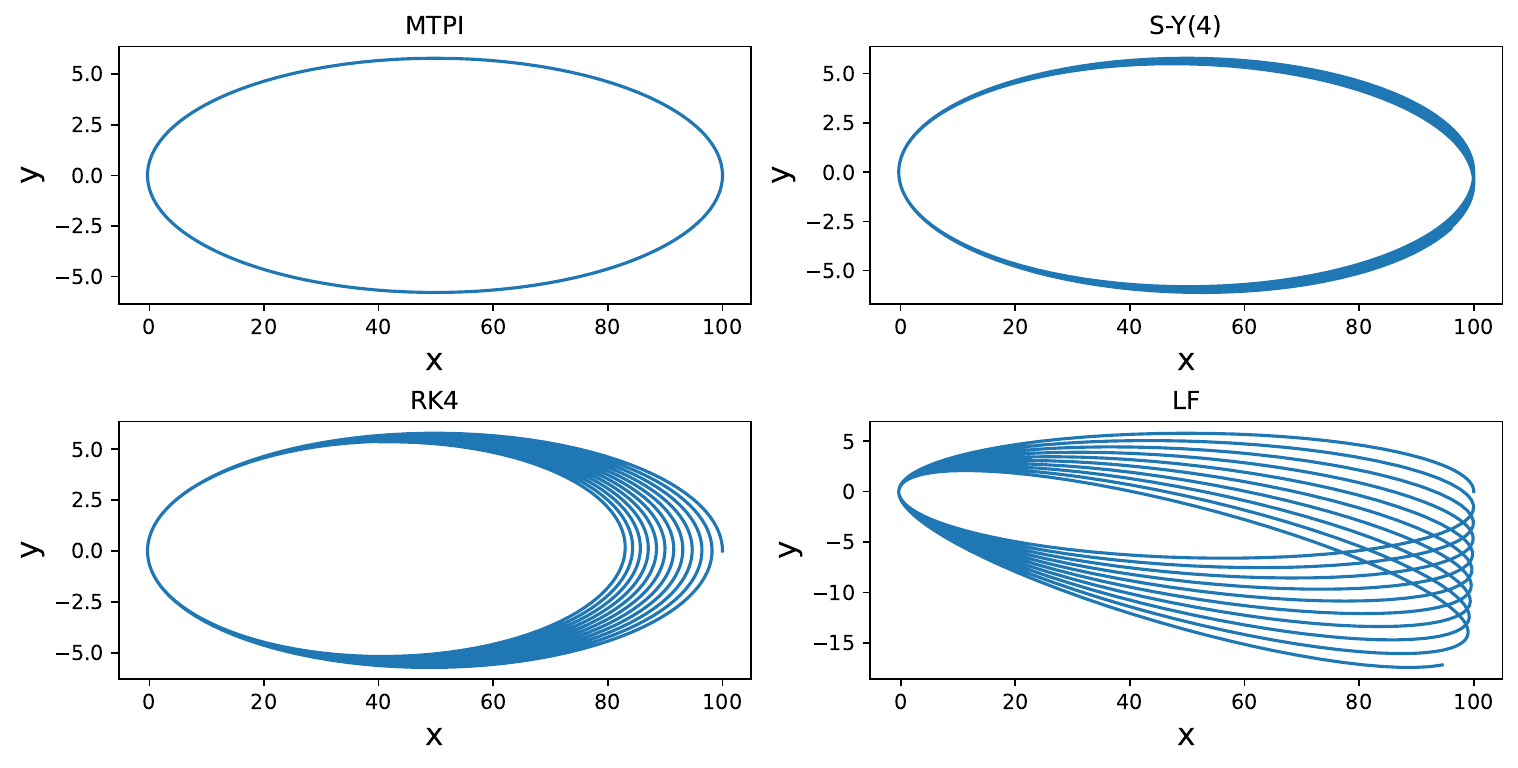} 
	\caption{Trajectories computed by the algorithms, projected onto the $x$–$y$ plane; TPI and MTPI yield almost identical orbits.}
	\label{fig:orbit}
\end{figure}

For the orbital shape, both constant angular-increment methods yield a stable ellipse. Other methods show unwanted numerical precession or orbit deformation, see Fig.~\ref{fig:orbit}.

\section{Conclusions}

The numerical scheme (MTPI) proposed in this Letter exhibits excellent accuracy and exactly preserves all first integrals (up to round-off errors), which in turn ensures exact preservation of orbital trajectories.
MTPI is fully explicit and therefore computationally efficient.
In our tests, it allowed time steps about 14 times larger than those used in the Runge–Kutta and Suzuki–Yoshida fourth-order schemes, yet still achieved accuracy superior by several orders of magnitude.
An additional remarkable feature of MTPI — shared also by the second constant–angular–increment integrator, TPI — is its ability to reproduce the exact time evolution along the trajectory.
Further studies are in progress to explore applications of this approach to many-body problems.

\appendix 

\section{Technical details} 

This appendix provides the derivations and proofs omitted from the main text, offering further insight into the principal results of the paper.

\subsection{Constant angular increment property}
\label{app-true}

We rewrite some of equations \rf{Kep-dis} in slightly modified form:
\be \ba{l} \dis  \label{Kep-dis-mod}
\pmb{p}_{n+1} = \pmb{p}_n - \lambda_n \rr_{n+1} \ ,  \\[2ex]\dis
h_n = \left( \displaystyle 2 \cos 2 \delta \ \frac{r_n}{r_{n+1}} - 1 + \frac{h_n \lambda_n}{ m } \right) h_{n+1} \ , \\[3ex]\dis
\frac{h_n \rr_{n+2}}{h_{n+1}} =  \left( 1 + \frac{h_n}{h_{n+1}} - \frac{\lambda_n h_n}{m} \right) \rr_{n+1} - \rr_n \ ,
\ea \ee
where
\be  \label{lambda_n}
\lambda_n = \frac{k h_n  }{   r_{n+1}^2 r_n \cos\delta} ,
\ee
Substituting $h_n/h_{n+1}$ from the second into the last equation of \rf{Kep-dis-mod} yields
\be  \label{r_n+2}
\rr_{n+2} = \frac{h_{n+1}}{h_n} \left( \frac{2 r_n \rr_{n+1} \cos 2\delta}{r_{n+1}} \ - \rr_n   \right) . 
\ee
We define 
\be  \label{delta_n}
\cos 2\delta_n := \frac{\rr_n \cdot \rr_{n+1}}{r_n r_{n+1}} .
\ee
We are going to show that $\delta_n$ does not depend on $n$ provided that we take into account $\delta_0 = \delta$  (which is a consequence of   \rf{delta_n} and \rf{delta}). Indeed, from \rf{r_n+2} we compute
\be  \label{delta_n+1}
\cos 2\delta_{n+1} = \frac{\rr_{n+1} \cdot \rr_{n+2}}{r_{n+1} r_{n+2}} = \frac{h_{n+1} r_n (2   \cos 2\delta \ -   \cos 2\delta_n)}{h_n  r_{n+2} } .
\ee
Then, again from \rf{r_n+2}, we get
\be  \label{r_n+2 r_n} 
r_{n+2} = \sqrt{\rr_{n+2} \cdot \rr_{n+2}}=  \frac{h_{n+1}}{h_n}  \, r_n
\ee
where we have used \rf{delta_n} to replace $\rr_{n+1}\cdot \rr_n$. Thus we can eliminate $r_{n+2}$ from \rf{delta_n+1} arriving at
\be
\cos 2\delta_{n+1} = 2   \cos 2\delta \ -   \cos 2\delta_n  
\ee
for  $n=0,1,2\ldots$.  Taking into account $\delta \equiv \delta_0$ we get, by mathematical induction,  $\delta_n = \delta$ for any $n$.  Therefore, the angle between $\pmb{q}_n$ and $\pmb{q}_{n+1}$ (which coincides with the angle between  $\rr_n$ and $\rr_{n+1}$) equals $2 \delta$ for all $n \in \mathbb N$.

\subsection{Preservation of first integrals}

\subsubsection{Angular momentum.} 

First, we express $\pmb{L}_n$ in terms of $\rr_n$ and $\rr_{n+1}$.

\be  \label{L-r-r}
\pmb{L}_n = \pmb{q}_n \times \pmb{p}_n = \frac{r_{n+1} \rr_n + r_n \rr_{n+1}}{r_n + r_{n+1}} \times  \frac{ m (\rr_{n+1} - \rr_n)}{h_n} =  \frac{m}{h_n} \ \rr_n \times \rr_{n+1} .
\ee
Then, taking into account  \rf{r_n+2}, we obtain  
\be
\pmb{L}_{n+1} = \frac{m}{h_{n+1}} \rr_{n+1} \times \rr_{n+2} = \frac{m}{h_{n+1}}  \frac{h_{n+1}}{h_n}  \rr_n \times \rr_{n+1}  = \pmb{L}_n 
\ee

\subsubsection{Energy.} 

Expressing $\pmb{p}_{n+1}$ and $\pmb{p}_n$ in terms of $\rr_n$ and $\rr_{n+1}$, we obtain
\be
\frac{\pmb{p}_{n+1}^2 - \pmb{p}_n^2}{2 m} = \frac{\lambda_n r_{n+1}^2}{2 h_n} \left( \frac{\lambda_n h_n}{m}  - 2 +  \frac{2 r_n  \cos 2\delta}{r_{n+1}}\right)  = \frac{\lambda_n r_{n+1}^2}{2 h_n} \left( \frac{h_n}{h_{n+1}}  - 1 \right) .
\ee
where we used \rf{Kep-dis-mod}. 
On the other hand, from
\be  \label{qn} 
\pmb{q}_n := \frac{ r_{n+1} \rr_n +  r_n \rr_{n+1}}{r_n + r_{n+1}} \quad \text{and} \quad \pmb{q}_n^2 = \frac{r_{n+1}^2  r_n^2 (2 + 2 \cos 2\delta )     }{(r_n + r_{n+1})^2}  
\ee
we obtain
\be  \label{1/q}
\frac{1}{q_n} = \frac{1}{2\cos\delta} \left( \frac{1}{r_n} + \frac{1}{r_{n+1}}  \right) .
\ee
Then, by \rf{r_n+2 r_n}, we have
\[
\frac{2\cos\delta}{q_n} - \frac{2\cos\delta}{q_{n+1}} \equiv  \frac{1}{r_n} - \frac{1}{r_{n+2}} = \frac{1}{r_n} \left( 1 - \frac{h_n}{h_{n+1}}  \right) .
\]
Therefore, from \rf{ene}, we obtain
\[
E_{n+1} - E_n = \left( \frac{h_n}{h_{n+1}}  - 1 \right) \left( \frac{\lambda_n r_{n+1}^2}{2 h_n}  - \frac{k}{2 r_n \cos\delta} \right) ,
\]
which equals zero for $\lambda_n$ given by \rf{lambda_n}.

\subsubsection{Laplace-Runge-Lenz vector.} 

Taking into account $\pmb{L}_{n+1}=\pmb{L}_n$ and the first equation of \rf{Kep-dis-mod}, we obtain from \rf{LRL}:
\be \label{shift A}
\pmb{A}_{n+1} - \pmb{A}_n = -  \frac{\lambda_n}{m} \ \rr_{n+1} \times \pmb{L}_n - k \left( \hat{\pmb{q}}_{n+1}  - \hat{\pmb{q}}_n           \right) ,
\ee
where the hat denotes a unit vector, i.e., $\hat{\pmb{q}}_n$  is  the unit vector in the direction of $\pmb{q}_n$. Then, by \rf{L-r-r},
\[
\frac{\lambda_n}{m} \ \rr_{n+1} \times \pmb{L}_n = 
\frac{\lambda_n}{h_n} \ (r_{n+1}^2  \rr_n - (\rr_n \cdot \rr_{n+1})  \rr_{n+1}) 
\]
Therefore, substituting into \rf{shift A}  $\lambda_n$ from \rf{lambda_n}, we have
\be  \label{A-A}
\frac{\pmb{A}_{n+1} - \pmb{A}_n}{k} =  \frac{\hat{\pmb{r}}_{n+1} \cos 2\delta - \hat{\pmb{r}}_n}{\cos\delta}  + \hat{\pmb{q}}_n  -  \hat{\pmb{q}}_{n+1} .
\ee
The righ-hand side can be easily shown to vanish. Indeed, from \rf{qn} and \rf{1/q} we obtain:
\be
\hat{\pmb{q}}_n = \frac{1}{2\cos\delta} \left( \hat{\rr}_{n+1} + \hat{\rr}_n \right) 
\ee
and, as a consequence,
\be  \label{qn-qn+2}
\hat{\pmb{q}}_n  -  \hat{\pmb{q}}_{n+1} = \frac{1}{2\cos\delta} \left( \hat{\rr}_n  -  \hat{\rr}_{n+2} \right) .
\ee
Then, dividing \rf{r_n+2} by \rf{r_n+2 r_n} we obtain
\be   \label{hat rn+2}
\hat{\rr}_{n+2}  =  2   \cos 2\delta\ \hat{\rr}_{n+1} - \hat{\rr}_n 
\ee
From \rf{qn-qn+2} and \rf{hat rn+2} it follows immediately that the right hand side of \rf{A-A} vanishes.

\subsection{First steps of the algorithm} 
\label{app-first}

The algorithm \rf{Kep-dis} needs $\pmb{q}_0$ and $\pmb{p}_0$ as initial data, while it apparently  starts rather from $\rr_0$ and $\pmb{p}_0$. We will derive equations \rf{Sn-rn}  which express $\rr_n$ in terms of   $\pmb{q}_n$ and $\pmb{p}_n$. First, we use \rf{qn} and the first equation of \rf{Kep-dis} to obtain
\be  \label{aaa}
\pmb{q}_n = \frac{ r_{n+1} \rr_n +  r_n \left( \rr_n + \frac{h_n \pmb{p}_n}{m} \right)}{r_n + r_{n+1}} = \rr_n + \frac{h_n r_n }{r_n + r_{n+1}} \ \pmb{p}_n .
\ee
Then, we will express the coefficient by $\pmb{p}_n$ in terms of $q_n$ and $S_n$. From \rf{Sn-rn} and \rf{Kep-dis} we obtain: 
\be  \label{qqSS} 
q_n S_n = \frac{h_n \pmb{q}_n \cdot \pmb{p}_n}{m} = 
\pmb{q}_n \cdot (\rr_{n+1} - \rr_n) = 
 \frac{2 (r_{n+1} - r_n) r_n r_{n+1} \cos^2 \delta}{ r_{n+1} + r_n} \ .  
\ee
Using  \rf{1/q} and \rf{qqSS}, and then once more \rf{1/q} itself,  we obtain
\be \label{r_n-r_n+1}
S_n  =  (r_{n+1} - r_n) \cos\delta \ , \qquad
r_{n+1} = \dis \frac{r_n q_n}{2 r_n \cos\delta - q_n} \ .
\ee  
Eliminating $r_{n+1}$ from the last two equations, we get the following quadratic equation
\be
r_n^2 \cos\delta + r_n \left( S_n  - q_n \right) - \frac{S_n  q_n}{ 2 \cos\delta}  = 0 \ ,
\ee
which has two solutions for $r_n$:
\be  \label{rn-pm}
r_n = \frac{q_n -  S_n \pm \sqrt{q_0^2 + S_n^2}}{2\cos\delta} \ .
\ee
\enlargethispage{2\baselineskip}
Using the first equation of \rf{r_n-r_n+1} we compute the corresponding $r_{n+1}$. Thus 
\be  \label{r_n+r_n+1}
r_{n+1} +r_n = \frac{q_n + \sqrt{q_0^2 + S_n^2}}{2\cos\delta} 
\ee
where we have selected the upper sign (plus) in the formula \rf{rn-pm}, because the right-hand side of \rf{r_n+r_n+1} has to be positive. Finally, using \rf{aaa}, \rf{r_n-r_n+1} and \rf{r_n+r_n+1}, we obtain the second equation of \rf{Sn-rn}.

%
%


\bibliographystyle{unsrtabbrv}

\bibliography{CJ-Kepler}

\end{document}